\documentclass[12pt, reqno]{amsart}
\usepackage{amsmath, amsthm, amscd, amsfonts, amssymb, graphicx, color, enumerate}
\usepackage[bookmarksnumbered, colorlinks, plainpages]{hyperref}
\usepackage[all,cmtip]{xy}
\hypersetup{colorlinks=true,linkcolor=red, anchorcolor=green, citecolor=cyan, urlcolor=red, filecolor=magenta, pdftoolbar=true}

\allowdisplaybreaks[4]

\newtheorem{theorem}{Theorem}[section]
\newtheorem{lemma}[theorem]{Lemma}

\newtheorem{corollary}[theorem]{Corollary}
\theoremstyle{definition}
\newtheorem{definition}[theorem]{Definition}
\newtheorem{example}[theorem]{Example}

\newtheorem{question}[theorem]{Question}

\numberwithin{equation}{section}
\newtheorem{claim}{Claim}

\newcommand{\Lip}{\mathcal L}
\newcommand{\Ran}{\operatorname{Ran}}
\newcommand{\codim}{\operatorname{codim}}
\newcommand{\ind}{\operatorname{ind}}
\begin{document}
\setcounter{page}{1}

\title[Fredholm weighted composition operators on the Lipschitz space]{Fredholm weighted composition operators on the Lipschitz space of a tree}

\author[J. Kong]{Junkai Kong}

\address{School of Mathematics and Statistics, Xidian University, Xi'an 710126, P.R. China.}

\email{\textcolor[rgb]{0.00,0.00,0.84}{xskongjunkai@163.com}}

\author[L. Liu]{Lei Liu$^{*}$}

\address{School of Mathematics and Statistics, Xidian University, Xi'an 710126, P.R. China.}

\email{\textcolor[rgb]{0.00,0.00,0.84}{leiliu@mail.xidian.edu.cn}}

\makeatletter
\@namedef{subjclassname@2020}{\textup{2020} Mathematics Subject Classification}
\makeatother


\date{Received: xxxxxx; Revised: yyyyyy; Accepted: zzzzzz.
\newline \indent $^{*}$ Corresponding author}

\begin{abstract}
Let $\mathcal{L}$ be the Lipschitz space of a tree. In this paper, we establish a complete characterization of the weighted composition operator $W_{\psi,\varphi}$ to be Fredholm on $\mathcal{L}$. Applying our characterization, we describe the connection between the Fredholmness of $W_{\psi,\varphi}$ and that of the composition operator $C_\varphi$.
\end{abstract} \maketitle

\section{Introduction}

Let $X$ be a Banach space consisting of complex-valued functions on a non-empty set $\Omega$. For a self-map $\varphi$ of $\Omega$ and a complex-valued function $\psi$ on $\Omega$, the corresponding weighted composition operator $W_{\psi,\varphi}$ is defined as
$$
W_{\psi,\varphi}(f)=\psi(f\circ\varphi) \qquad \forall f\in X.
$$
If $\psi\equiv1$, then $W_{\psi,\varphi}$ reduces to a composition operator $C_\varphi$ defined by $C_\varphi(f)=f\circ\varphi$. If $\varphi$ is the identity map, then $W_{\psi,\varphi}$ reduces to a multiplication operator $M_\psi$ defined by $M_\psi(f)=\psi f$.

Recall that a bounded linear operator $A$ on a Banach space is called a Fredholm operator if $A$ has closed range, $\dim(\ker A)<\infty$, and $\codim(\Ran A)<\infty$. It is trivial that every invertible operator on a Banach space is Fredholm. Fredholmness of (weighted) composition operators has been extensively investigated on Hardy spaces, weighted Bergman spaces, Hardy-Sobolev spaces and other general Banach spaces of holomorphic functions (for example, refer to \cite{bo,chl,chz,he,hlns,mz} and the references therein).

In recent years, there have been a number of papers on the study of (weighted) composition operators on function spaces over discrete structures such as trees. Allen, Colonna, and Easley \cite{ace} characterized the boundedness, compactness, isometries and spectrum of the composition operators on the Lipschitz space of a tree. Colonna and  Mart\'{i}nez-Avenda$\tilde{\mathrm{n}}$o \cite{cma} characterized the boundedness, spectrum and hypercyclicity of composition operators on the little Lipschitz space of a tree. In \cite{msk} the bounded and compact weighted composition operators between discrete weighted Hardy spaces on rooted trees were investigated. However, to the best of our knowledge, little is known about Fredholm (weighted) composition operators on function spaces of a tree. This naturally leads to the following question.
\begin{question}\label{que}
Is it possible to establish a characterization of $W_{\psi,\varphi}$ to be a Fredholm operator on the Lipschitz space of a tree?
\end{question}

In this paper, we give a confirmative answer to the above question.

To state our main result, we first introduce some notation. Let $T$ be a tree, and let $\Lip$ denote the Lipschitz space of $T$. Let $\varphi:T\to T$ and $\psi:T\to\mathbb C$. We set
$$
R=\varphi(T).
$$
For every $w\in R$, define
$$
F_w=\varphi^{-1}(w),\qquad N_w=|F_w|.
$$
We also set
$$
N=|T\setminus R|, \qquad K=\sum_{w\in R}(N_w-1)
$$
and
$$
Z_\psi=\{x\in T:\psi(x)=0\}.
$$
Define the set of ``good vertices" by
$$
G=\{w\in R:F_w=\{s(w)\},\ \psi(s(w))\neq0\},
$$
where $s(w)$ is the unique preimage	of $w$ under $\varphi$.
For any $u,v\in G$ with $u\sim v$, and
$
x=s(u), y=s(v), z=x\wedge y,
$
we define
$$
\Lambda_{\psi,\varphi}(u,v)
=(|z|+1)\left|\frac1{\psi(x)}-\frac1{\psi(y)}\right|
+\frac{d(z,x)}{|\psi(x)|}
+\frac{d(z,y)}{|\psi(y)|}.
$$

We are now ready to state our main result, Theorem \ref{thm} below, which answers Question \ref{que}.

\begin{theorem}\label{thm}
Let $W_{\psi,\varphi}$ be a bounded weighted composition operator on $\Lip$. Then $W_{\psi,\varphi}$ is Fredholm if and only if the following conditions hold:
\begin{enumerate}[(i)]
\item $N<\infty$;
\item $|Z_\psi|<\infty$;			
\item for any $w\in R$, $N_w<\infty$ and $K<\infty$;
\item $\sup\{\Lambda_{\psi,\varphi}(u,v): u,v\in G \ \text{with} \ \ u\sim v\}<\infty$.
\end{enumerate}
In this case, $\ind W_{\psi,\varphi}=N-K$.
\end{theorem}

The exposition of the paper is organized as follows. We refer to  Section 2 for some preliminaries on trees and the definitions of the Lipschitz space of a tree. Section 3 is devoted to the proof of Theorem ~\ref{thm}. In Section 4 we characterize invertible weighted composition operators and Fredholm composition operators as an application of Theorem~\ref{thm}. We then give a characterization of the relationship between the Fredholmness of the weighted composition operator $W_{\psi,\varphi}$ and that of the composition operator $C_\varphi$.

\section{preliminaries}

A tree is a locally finite connected and simply connected graph $T$, which, as a set, we identify with the collection of its vertices. Two vertices $u$ and $v$ are said to be neighbors if there is an edge connecting them, and we adopt the notation $u \sim v$. The degree of a vertex is the number of its neighbors. A path between two vertices $u$ and $w$ is a finite sequence of distinct vertices $[u, u_1, u_2, \ldots, u_{k}, w]$ such that $u \sim u_1 \sim u_2 \sim \cdots \sim u_{k} \sim w$. Since $T$ is a tree, the path between two vertices is unique.
	
Given a tree $T$ rooted at $o$ and a vertex $u\in T$, a vertex $v$, with $v \neq u$, is called a descendant of $u$ if $u$ lies in the unique path from $o$ to $v$. The vertex $u$ is then called an ancestor of $v$. We denote by $v^-$ the neighbor of $v$ that is an ancestor of $v$, and call it the parent of $v$. The vertex $v$ is called a child of $v^-$. The length of a path between two vertices $u$ and $w$ is the number of edges connecting them. The distance $d(u, w)$ between vertices $u$ and $w$ is the length of the path connecting $u$ to $w$. In particular, the length of a path joining the root $o$ to a vertex $w$ is called the length of $w$ and is denoted by $|w|$. A function on a tree is a complex-valued function on its vertex set.
	
Throughout this paper, we assume that $T$ is a rooted tree with root $o$ whose vertices, save possibly $o$, have degree greater than 1, and hence the tree is infinite. $[u,w]$ denotes the path $u \sim u_1 \sim u_2 \sim \cdots \sim u_{k} \sim w$, while $(u,w]$ denotes the set $[u,w]\backslash\{u\}$. We adopt the notation $T^*$ for the set $T \setminus \{o\}$ and $x\wedge y$ for the greatest common ancestor of $x$ and $y$.
	
\begin{definition}
The Lipschitz space of a tree $T$ is the set $\mathcal{L}$ of functions $f: T \to \mathbb{C}$ such that
$$
\beta_f:=\sup_{\substack{v \neq w \\ v, w \in T}} \frac{|f(v) - f(w)|}{d(v, w)} < \infty.
$$
\end{definition}

It can be shown in \cite{ce} that
$$
\sup_{v \in T^*} |f(v) - f(v^-)| = \sup_{\substack{v \neq w \\ v, w \in T}} \frac{|f(v) - f(w)|}{d(v, w)}.
$$
Thus we can also define the Lipschitz space of a tree $T$ to be the set $\mathcal{L}$ of functions $f: T \to \mathbb{C}$ satisfying
$$
\sup_{v \in T^*} |f(v) - f(v^-)| < \infty.
$$
Colonna and Easley \cite{ce} showed that $\mathcal{L}$ is a Banach space when endowed with the norm
$$
\|f\| := |f(o)|+\sup_{v \in T^*} |f(v) - f(v^-)|.
$$
In this paper, we endow the Lipschitz space $\mathcal{L}$ with the equivalent norm as follows.
$$
\|f\|_{\mathcal{L}} := \max\left\{ |f(o)|, \sup_{v \in T^*} |f(v) - f(v^-)| \right\}.
$$
For any $f\in\Lip$, define
$$
D:\Lip\longrightarrow\ell^\infty(T)
$$
by
$$
(Df)(v) = \begin{cases}
f(v) - f(v^-), & \text{if } v \neq o, \\
f(o), & \text{if } v = o.
\end{cases}
$$
It follows that
$$
\|f\|_{\mathcal{L}} = \sup_{v \in T} |(Df)(v)|.
$$
This yields the mapping $D$ is an isometric isomorphism. Its inverse is given by
$$
f(w)=\sum_{z\in[o,w]}Df(z),\qquad w\in T.
$$
Consequently,
\begin{equation}\label{eq.2.1}
|f(w)|\le (|w|+1)\|f\|_{\Lip},\qquad w\in T.
\end{equation}
This implies that every evaluation functional $\delta_w(f)=f(w)$ belongs to $\Lip^*$.
	
We conclude this section with a useful lemma that will be used repeatedly.

\begin{lemma}\label{lem.2.2}
Let $x,y\in T$, $a,b\in\mathbb C$ and $z=x\wedge y$. Define
$$
\Phi_{a,b;x,y}(f)=af(x)-bf(y),\qquad f\in\Lip.
$$
Then
$$
\|\Phi_{a,b;x,y}\|=(|z|+1)|a-b|+d(z,x)|a|+d(z,y)|b|.
$$
\end{lemma}

\begin{proof}
For any $f\in\Lip(T)$, we have
$$
f(x)=\sum_{v\in[o,z]}Df(v)+\sum_{v\in(z,x]}Df(v)
$$
and
$$
f(y)=\sum_{v\in[o,z]}Df(v)+\sum_{v\in(z,y]}Df(v).
$$
Hence
\begin{equation}\label{eq.2.2}
\Phi_{a,b;x,y}(f)=(a-b)\sum_{v\in[o,z]}Df(v)+a\sum_{v\in(z,x]}Df(v)
-b\sum_{v\in(z,y]}Df(v).
\end{equation}
Since $D$ is an isometric isomorphism onto $\ell^\infty(T)$, it follows from \eqref{eq.2.2} that
\begin{equation}\label{eq.2.3}
\|\Phi_{a,b;x,y}\|\le(|z|+1)|a-b|+d(z,x)|a|+d(z,y)|b|.
\end{equation}

Conversely, we may define $Df$ separately on the three disjoint sets $[o,z]$, $(z,x]$ and $(z,y]$ by setting $|Df(v)|=1$ at every vertex where the summands in \eqref{eq.2.2} have nonzero coefficients, and choosing the argument of $Df(v)$ so that all terms $(a-b)Df(v)$ for $v\in[o,z]$, $aDf(v)$ for $v\in(z,x]$, and $-bDf(v)$ for $v\in(z,y]$ have the same argument. Let $Df(v)=0$ elsewhere.
Hence $\|Df\|_\infty=1$, which implies $\|f\|_{\mathcal L}=1$ by the isometry of $D$. Therefore,
$$
\big|\Phi_{a,b;x,y}(f)\big|=\big(|z|+1\big)|a-b|+d(z,x)|a|+d(z,y)|b|,
$$
and consequently
\begin{equation}\label{eq.2.4}
\|\Phi_{a,b;x,y}\| \ge \big(|z|+1\big)|a-b| + d(z,x)|a| + d(z,y)|b|.
\end{equation}
Combining \eqref{eq.2.3} and \eqref{eq.2.4}, we obtain equality. This completes the proof.
\end{proof}

\section{Proof of Theorem \ref{thm}}

We need a lemma to prove Theorem \ref{thm}.

Let $u,v\in G$ with $u\sim v$.
Define a linear functional by
$$
L_{u,v}(h)=\frac{h(x)}{\psi(x)}-\frac{h(y)}{\psi(y)},
\qquad h\in\Lip,
$$
where $x=s(u)$ and $y=s(v)$. By Lemma~\ref{lem.2.2}, we have $L_{u,v}\in\Lip^*$.
	
\begin{lemma}\label{lem.3.1}
Let $\{u_n\sim v_n\}_{n=1}^\infty\subseteq G$ be a sequence of neighbors satisfying $\min\{|u_n|,|v_n|\}\to\infty$ as $n\to\infty$.
If
$$
\|L_{u_n,v_n}\|\to\infty
$$
as $n\to\infty$, then the sequence
$$
\big\{\frac{L_{u_n,v_n}}{\|L_{u_n,v_n}\|}\big\}_{n=1}^\infty
$$
has no norm convergent subsequence in $\Lip^*$.
\end{lemma}

\begin{proof}
Write
$$
x_n=s(u_n),\qquad y_n=s(v_n),\qquad z_n=x_n\wedge y_n,
$$
and
$$
a_n=\frac1{\psi(x_n)},\qquad b_n=\frac1{\psi(y_n)}.
$$
Under the isometric isomorphism $D:\Lip(T)\to\ell^\infty(T)$, the linear functional $L_{u_n,v_n}$ is represented by $l_n\in\ell^1(T)$, where
$$
l_n(w)=
\begin{cases}
a_n-b_n, & w\in[o,z_n],\\
a_n, & w\in(z_n,x_n],\\
-b_n, & w\in(z_n,y_n],\\
0, & \text{otherwise}.
\end{cases}
$$
Thus, by Lemma~\ref{lem.2.2}, we have
\begin{equation}\label{eq.3.1}
\|l_n\|_1=(|z_n|+1)|a_n-b_n|+d(z_n,x_n)|a_n|+d(z_n,y_n)|b_n|
=\|L_{u_n,v_n}\|.
\end{equation}
Set
$$
\delta_n=\frac{l_n}{\|l_n\|_1}.
$$
Note that the canonical embedding $J:\ell^1(T)\longrightarrow(\ell^\infty(T))^*$
is isometric. It suffices to show that $\{\delta_n\}_{n=1}^{\infty}$ has no norm convergent subsequence in $\ell^1(T)$.
		
We proceed by contradiction. Assume that there exists subsequence $\{\delta_{n_k}\}_{k=1}^{\infty}$ converging in $\ell^1(T)$. Without loss of generality, we may further assume that either $|z_{n_k}|\to\infty$ or $\{|z_{n_k}|\}_{k=1}^{\infty}$ is bounded.
		
For the case $|z_{n_k}|\to\infty$, fix $M\ge0$ and let
$$
B_M=\{z\in T:|z|\le M\}.
$$
For sufficiently large $k$, the sets $(z_{n_k},x_{n_k}]$ and $(z_{n_k},y_{n_k}]$ are disjoint from $B_M$. If $a_{n_k}=b_{n_k}$, then
$$
\sum_{w\in B_M}|\delta_{n_k}(w)|=0.
$$
If $a_{n_k}\neq b_{n_k}$, then by \eqref{eq.3.1}, we have
$$
\sum_{w\in B_M}|\delta_{n_k}(w)|\le\frac{(M+1)|a_{n_k}-b_{n_k}|}{(|z_{n_k}|+1)|a_{n_k}-b_{n_k}|}
=\frac{M+1}{|z_{n_k}|+1},
$$
which implies		
$$\label{eq:mass-escape-case1}
\sum_{w\in B_M}|\delta_{n_k}(w)|\to0
$$
as $k\to\infty$.
		
For the case $\{|z_{n_k}|\}_{k=1}^{\infty}$ is bounded, we claim
\begin{equation}\label{eq.3.2}
|x_{n_k}|\to\infty,
\qquad
|y_{n_k}|\longrightarrow\infty
\end{equation}
as $k\to\infty$. Indeed, if $B\subseteq T$ is finite, then
$$
\{w\in G:s(w)\in B\}\subseteq\varphi(B),
$$
and the latter set is finite. Since $\min\{|u_{n_k}|,|v_{n_k}|\}\to\infty$ as $k\to\infty$, \eqref{eq.3.2} follows.
Set
$$
r_k=d(z_{n_k},x_{n_k}),\qquad t_k=d(z_{n_k},y_{n_k}).
$$
The boundedness of $\{|z_{n_k}|\}_{k=1}^{\infty}$ and \eqref{eq.3.2} imply
$$
r_k\to\infty,
\qquad
t_k\to\infty
$$
as $k\to\infty$.
On the one hand, for each fixed $M\ge0$, there exists $C_M>0$, independent of $k$, such that
$$
\sum_{w\in B_M}|l_{n_k}(w)|\le C_M(|a_{n_k}|+|b_{n_k}|).
$$
On the other hand, it follows from \eqref{eq.3.1} that
$$
\|l_{n_k}\|_1\ge r_k|a_{n_k}|+t_k|b_{n_k}|\ge\min\{r_k,t_k\}(|a_{n_k}|+|b_{n_k}|).
$$
Therefore
$$\label{eq:mass-escape-case2}
\sum_{w\in B_M}|\delta_{n_k}(w)|\le\frac{C_M}{\min\{r_k,t_k\}}\to0
$$
as $k\to\infty$.
		
Thus, in both cases and for a fixed $M\ge0$, we obtain that
\begin{equation}\label{eq.3.3}
\sum_{w\in B_M}|\delta_{n_k}(w)|\to0
\end{equation}
as $k\to\infty$.
Moreover, if $\delta_{n_k}\to\delta$ in $\ell^1(T)$, then for any $M\ge0$, we get
$$
\left|
\sum_{w\in B_M}|\delta_{n_k}(w)|-\sum_{w\in B_M}|\delta(w)|
\right|
\le
\|\delta_{n_k}-\delta\|_1\to0
$$
as $k\to\infty$. Combined with \eqref{eq.3.3}, this yields
$$
\sum_{w\in B_M}|\delta(w)|=0
$$
for all $M\ge0$. Since $T=\bigcup_{M\ge0}B_M$, it follows that $\delta=0$. This contradicts
$$
\|\delta_{n_k}\|_1=1
$$
for all $k$.
Hence $\{\delta_n\}_{n=1}^{\infty}$ has no norm convergent subsequence in $\ell^1(T)$.
\end{proof}

Now we prove Theorem \ref{thm}.

\begin{proof}[Proof of Theorem \ref{thm}]
``$\Rightarrow$" Suppose that $N=\infty$. Then, for every $w\in T\setminus R$, the characteristic function $\chi_{\{w\}}$ belongs to $\Lip$ and satisfies
$$
W_{\psi,\varphi}\chi_{\{w\}}=0.
$$
It is clear that the family $\{\chi_{\{w\}}:w\in T\setminus R\}$ is linearly independent. This yields $\ker W_{\psi,\varphi}$ is infinite dimensional, contradicting Fredholmness. Hence $N<\infty$.
		
Moreover, set $Y=\Ran W_{\psi,\varphi}$. If $x\in Z_\psi$, then
$$
(W_{\psi,\varphi}f)(x)=0, \qquad f\in\Lip,
$$
which implies the evaluation function $\delta_x\in Y^\perp$. The family of evaluation functions at distinct vertices is linearly independent. Indeed, if $x_1,\dots,x_m$ are distinct and $\sum_{j=1}^m\alpha_j\delta_{x_j}=0$, then evaluating against the characteristic function $\chi_{\{x_j\}}$ gives $\alpha_j=0$. Since $\codim(Y)<\infty$, it follows that $\dim(Y^\perp)=\codim(Y)$, and hence $|Z_\psi|<\infty$.

Furthermore, assume that there exists $w\in R$ such that $|N_w|=\infty$. Since $Z_\psi$ is finite, there exists $b_w\in F_w$ with $\psi(b_w)\neq0$. For any $x\in F_w\setminus\{b_w\}$, define
\begin{equation}\label{eq.3.4}
\lambda_{w,x}(h)=h(x)-\frac{\psi(x)}{\psi(b_w)}h(b_w), \qquad h\in\Lip.
\end{equation}
If $h=W_{\psi,\varphi}f$, then we have $\lambda_{w,x}(h)=0$, which implies $\lambda_{w,x}\in Y^\perp$. The family
$$
\{\lambda_{w,x}:x\in F_w\setminus\{b_w\}\}
$$
is linearly independent. Indeed, suppose $x_1,\dots,x_m\in F_w\setminus\{b_w\}$ are distinct and $\sum_{j=1}^m\alpha_j\lambda_{w,x_j}=0$. Evaluating both sides against the characteristic function $\chi_{\{x_i\}}$ yields $\alpha_i=0$. This contradicts $\dim(Y^\perp)<\infty$. Hence $|N_w|<\infty$ for all $w\in R$.

Next suppose that $K=\infty$. Put
$$
Z_0=\{w\in R:\psi|_{F_w}\equiv0\}.
$$
Since
$$
\bigcup_{w\in Z_0}F_w\subset Z_\psi
$$
and $|Z_\psi|<\infty$, it follows that
$$
\sum_{w\in Z_0}(N_w-1)
$$
is finite. Moreover, for any $w\in R\setminus Z_0$ with $N_w>1$, take $b_w\in F_w$ such that $\psi(b_w)\neq0$ and define $\lambda_{w,x}$ by \eqref{eq.3.4} for $x\in F_w\setminus\{b_w\}$. We have proven that all these functionals belong to $Y^\perp$. Notice that $\{\lambda_{w,x}\}$ are jointly linearly independent for distinct $w$, since the sets $F_w$ are pairwise disjoint. If
$$
\sum_{w\in R\setminus Z_0}(N_w-1)=\infty,
$$
then this would contradict $\dim(Y^\perp)<\infty$. Therefore $K<\infty$.

Finally, since $W_{\psi,\varphi}$ is Fredholm, $Y$ is closed. Define
$$
\widetilde W_{\psi,\varphi}:
\Lip(T)/\ker W_{\psi,\varphi}\to Y
\qquad \text{by} \qquad
\widetilde W_{\psi,\varphi}(f+\ker W_{\psi,\varphi})=W_{\psi,\varphi}f,
$$
which is a bounded bijective linear mapping between Banach spaces. By the inverse mapping theorem, there exists $C>0$ such that
\begin{equation}\label{eq.3.5}
d(f,\ker W_{\psi,\varphi})\le C\|W_{\psi,\varphi}f\|_{\Lip}, \qquad f\in\Lip.
\end{equation}
Let $u,v\in G$ with $u\sim v$. If $k\in\ker W_{\psi,\varphi}$, then $k(u)=k(v)=0$. Indeed, for $x=s(u)$, we have
$$
0=(W_{\psi,\varphi}k)(x)=\psi(x)k(u).
$$
Since $\psi(x)\neq0$, it follows that $k(u)=0$. Using the same argument to $v$, we have $k(v)=0$.
Notice that if $h=W_{\psi,\varphi}f$, then
$$
L_{u,v}(h)=f(u)-f(v).
$$
Hence we arrive at
$$
|L_{u,v}(h)|=|(f-k)(u)-(f-k)(v)|\le\|f-k\|_{\Lip}
$$
for all $h\in Y$. Taking the infimum over $k\in\ker W_{\psi,\varphi}$ and using \eqref{eq.3.5}, we see that
\begin{equation}\label{eq.3.6}
\|L_{u,v}|_Y\|\le C.
\end{equation}
By the following isometric isomorphism
$$
\Lip(T)^*/Y^\perp\cong Y^*,
$$
we have
\begin{equation}\label{eq.3.7}
d(L,Y^\perp)=\|L|_Y\|,
\qquad L\in\Lip^*.
\end{equation}
Combining \eqref{eq.3.6} and \eqref{eq.3.7}, we obtain that
\begin{equation}\label{eq.3.8}
d(L_{u,v},Y^\perp)\le C
\end{equation}
for all $u, v\in G$ with $u\sim v$.

Assume that
$\sup\{\Lambda_{\psi,\varphi}(u,v): u,v\in G, u\sim v\}=\infty$.
By Lemma~\ref{lem.2.2}, we have $\Lambda_{\psi,\varphi}(u,v)=\|L_{u,v}\|$.
Since $T$ is locally finite, we may select sequences $u_n, v_n\in G$ satisfying $u_n\sim v_n$ and $\min\{|u_n|,|v_n|\}\to\infty$ such that
$$
\|L_{u_n,v_n}\|\to\infty
$$
as $n\to\infty$. By \eqref{eq.3.8}, there exists $q_n\in Y^\perp$ such that
$$
\|L_{u_n,v_n}-q_n\|\le C+1.
$$
Then
\begin{equation}\label{eq.3.9}
\left\|
\frac{L_{u_n,v_n}}{\|L_{u_n,v_n}\|}-\frac{q_n}{\|L_{u_n,v_n}\|}
\right\|
\to0
\end{equation}
as $n\to\infty$,
and
$$
\left\|\frac{q_n}{\|L_{u_n,v_n}\|}\right\|\le1+\frac{C+1}{\|L_{u_n,v_n}\|}.
$$
It follows that
$\big\{\frac{q_n}{\|L_{u_n,v_n}\|}\big\}_{n=1}^\infty$ is bounded in the finite dimensional space $Y^\perp$. Therefore, it has a norm convergent subsequence. By \eqref{eq.3.9}, the corresponding subsequence of
$\big\{\frac{L_{u_n,v_n}}{\|L_{u_n,v_n}\|}\big\}_{n=1}^\infty$ converges in norm, contradicting Lemma~\ref{lem.3.1}. Therefore, $\sup\{\Lambda_{\psi,\varphi}(u,v): u,v\in G, u\sim v\}<\infty$.

``$\Leftarrow$" Since $\bigcup_{w\in Z_0}F_w\subseteq Z_\psi$ and condition (ii), the set $Z_0$ is finite. Moreover, we have
$$\label{eq:G-complement}
T\setminus G=(T\setminus R)\cup\{w\in R:N_w>1\}\cup Z_0.
$$
By conditions (i)-(iii), we arrive at
\begin{equation}\label{eq.3.10}
|T\setminus G|<\infty.
\end{equation}

Next, define
$$
\mathcal E_{\psi,\varphi}=\left\{h\in\Lip(T):
h|_{F_w}\in\operatorname{span}\{\psi|_{F_w}\}\text{ for every }w\in R\right\}.
$$
For any $w\in R\setminus Z_0$, take $b_w\in F_w$ such that $\psi(b_w)\neq0$.
The condition on $F_w$ is equivalent to
\begin{equation}\label{eq.3.11}
h(x)-\frac{\psi(x)}{\psi(b_w)}h(b_w)=0,
\qquad x\in F_w\setminus\{b_w\}.
\end{equation}
For any $w\in Z_0$, the condition on $F_w$ is equivalent to
\begin{equation}\label{eq.3.12}
h(x)=0,
\qquad x\in F_w.
\end{equation}
It follows that each is the kernel condition of a continuous linear functional on $\Lip(T)$. Hence $\mathcal E_{\psi,\varphi}$ is closed.

Moreover, define
\begin{equation}\label{eq.3.13}
(\mathcal R_{\psi,\varphi}h)(w)=
\begin{cases}
\dfrac{h(b_w)}{\psi(b_w)}, & w\in R\setminus Z_0,\\[2mm]
0, & w\in(T\setminus R)\cup Z_0,
\end{cases}
\qquad h\in\mathcal E_{\psi,\varphi}.
\end{equation}
We show that $\mathcal R_{\psi,\varphi}:\mathcal E_{\psi,\varphi}\to\Lip$ is bounded.
Indeed, by condition (iv), we have
$$
M:=\sup_{\substack{u\sim v\\u,v\in G}}\Lambda_{\psi,\varphi}(u,v)<\infty.
$$
If $u,v\in G$ with $u\sim v$, then $b_u=s(u)$ and $b_v=s(v)$, so by \eqref{eq.3.13} and Lemma~\ref{lem.2.2}, we have
\begin{equation}\label{eq.3.14}
|(\mathcal R_{\psi,\varphi}h)(u)-(\mathcal R_{\psi,\varphi}h)(v)|=|L_{u,v}(h)|
\le M\|h\|_{\Lip}.
\end{equation}
Set
$$
E=\{\{u,v\}:u\sim v,\ \{u,v\}\cap(T\setminus G)\neq\varnothing\}.
$$
By \eqref{eq.3.10} and the local finiteness of $T$, the set $E$ is finite.
For any $w\in T$, define
$$
A_w=
\begin{cases}
\dfrac{|b_w|+1}{|\psi(b_w)|}, & w\in R\setminus Z_0,\\[2mm]
0, & w\in(T\setminus R)\cup Z_0.
\end{cases}
$$
Thus, by \eqref{eq.2.1} and \eqref{eq.3.13}, we have
\begin{equation}\label{eq.3.15}
|(\mathcal R_{\psi,\varphi}h)(w)|
\le A_w\|h\|_{\Lip},
\qquad w\in T.
\end{equation}
Consequently, for any $\{u,v\}\in E$, it follows that
\begin{equation}\label{eq.3.16}
|(\mathcal R_{\psi,\varphi}h)(u)-(\mathcal R_{\psi,\varphi}h)(v)|
\le(A_u+A_v)\|h\|_{\Lip}.
\end{equation}
Since $E$ is finite, we have
$$
M_{E}:=\max_{\{u,v\}\in E}(A_u+A_v)<\infty,
$$
with the convention $M_{E}=0$ if $E=\varnothing$. Also, it follows from \eqref{eq.3.15} that
$$
|(\mathcal R_{\psi,\varphi}h)(o)|\le A_o\|h\|_{\Lip}.
$$
Combining the above inequality with \eqref{eq.3.14} and \eqref{eq.3.16}, we obtain that
$$\label{eq:R-bounded}
\|\mathcal R_{\psi,\varphi}h\|_{\Lip}\le\max\{A_o,M,M_{E}\}\|h\|_{\Lip}.
$$
Therefore, $\mathcal R_{\psi,\varphi}:\mathcal E_{\psi,\varphi}\to\Lip$ is bounded.

Let $h\in\mathcal E_{\psi,\varphi}$ and $x\in T$, and put $w=\varphi(x)$. If $w\in R\setminus Z_0$, then by \eqref{eq.3.11} we have
$$
h(x)=\psi(x)\frac{h(b_w)}{\psi(b_w)}=\psi(x)(\mathcal R_{\psi,\varphi}h)(w).
$$
If $w\in Z_0$, then $h(x)=0$ by \eqref{eq.3.12} and $\psi(x)=0$. It follows that
$$
W_{\psi,\varphi}\mathcal R_{\psi,\varphi}h=h,
$$
which implies
$$
\mathcal E_{\psi,\varphi}\subseteq\Ran W_{\psi,\varphi}.
$$
Conversely, if $h=W_{\psi,\varphi}f$ and $x\in F_w$, then
$$
h(x)=\psi(x)f(w).
$$
This yields $h|_{F_w}\in\operatorname{span}\{\psi|_{F_w}\}$ for all $w\in R$, and hence $\Ran W_{\psi,\varphi}\subseteq\mathcal E_{\psi,\varphi}$. Therefore,
\begin{equation}\label{eq.3.17}
\Ran W_{\psi,\varphi}=\mathcal E_{\psi,\varphi}.
\end{equation}
In particular, we see that $\Ran W_{\psi,\varphi}$ is closed.

Now we compute $\ker W_{\psi,\varphi}$. For every $f\in\Lip$, $W_{\psi,\varphi}f=0$ if and only if $f(w)=0$ for all $w\in R\setminus Z_0$. Indeed, for any $w\in R\setminus Z_0$, we may choose $x\in F_w$ such that $\psi(x)\neq0$. Then
$$
0=(W_{\psi,\varphi}f)(x)=\psi(x)f(w),
$$
which implies $f(w)=0$. The converse is obvious.
It follows that $f\in \ker W_{\psi,\varphi}$ if and only if $f\in \operatorname{span}\{\chi_{\{w\}}: w\in(T\setminus R)\cup Z_0\}$.
This yields that
\begin{equation}\label{eq.3.18}
\dim(\ker W_{\psi,\varphi})=N+|Z_0|.
\end{equation}

Finally, we compute the codimension of $\Ran W_{\psi,\varphi}$. Let
$$
S=\bigcup_{\substack{w\in R\\N_w>1}}F_w\cup\bigcup_{w\in Z_0}F_w.
$$
By conditions (ii)-(iii), the set $S$ is finite. It is clear that the mapping
$$
\rho_S:\Lip(T)\longrightarrow\mathbb C^S, \qquad\rho_S(h)=h|_S,
$$
is surjective. For each $w$ with $F_w\subseteq S$, set
$$
V_w=
\begin{cases}
\operatorname{span}\{\psi|_{F_w}\}, & w\notin Z_0,\\
0, & w\in Z_0,
\end{cases}
$$
and
$$
V=\prod_{\substack{w\in R\\F_w\subseteq S}}V_w\subseteq\mathbb C^S.
$$
Then, by the definition of $\mathcal E_{\psi,\varphi}$, we have
$$
\mathcal E_{\psi,\varphi}=\rho_S^{-1}(V).
$$
Since $\rho_S$ is surjective, the following mapping
$$
\Lip(T)/\mathcal E_{\psi,\varphi}\to\mathbb C^S/V
$$
is an isomorphism. Hence, by \eqref{eq.3.17}, it follows that
\begin{align}
\codim(\Ran W_{\psi,\varphi})
&=\codim(\mathcal E_{\psi,\varphi})\notag\\
&=\codim (V)\notag\\
&=\sum_{w\in R\setminus Z_0}(N_w-1)+\sum_{w\in Z_0}N_w\notag\\
&=K+|Z_0|.
\label{eq.3.19}
\end{align}
Combining \eqref{eq.3.18} and \eqref{eq.3.19}, we obtain that $W_{\psi,\varphi}$ is Fredholm. In particular,
\begin{align*}
\ind W_{\psi,\varphi}&=\dim(\ker W_{\psi,\varphi})-\codim(\Ran W_{\psi,\varphi})\\
&=(N+|Z_0|)-(K+|Z_0|)\\
&=N-K.
\end{align*}
The proof is complete.
\end{proof}

\section{Applications}

In this section, as an application of Theorem~\ref{thm}, we first characterize invertible weighted composition operators and Fredholm composition operators on $\Lip$. Using our main result, we also investigate the relationship between the Fredholmness of weighted composition operators and that of composition operators.

We begin with a characterization of invertibility for weighted composition operators on $\Lip$.

\begin{theorem}\label{thm.4.1}
Let $W_{\psi,\varphi}$ be a bounded weighted composition operator
on $\Lip$. Then $W_{\psi,\varphi}$ is invertible if and only if
the following conditions hold:
\begin{enumerate}[(i)]
\item $\varphi:T\to T$ is a bijection;
\item $Z_\psi=\varnothing$;
\item $\sup\{\Lambda_{\psi,\varphi}(u,v): u,v\in T \ \text{with} \ \ u\sim v\}<\infty$.
\end{enumerate}
\end{theorem}

\begin{proof}
``$\Rightarrow$" Suppose $\psi(x)=0$ at some $x\in T$, then every function in
$\Ran W_{\psi,\varphi}$ vanishes at $x$, which contradicts surjectivity. Thus $Z_\psi=\varnothing$, and hence $Z_0=\varnothing$. Moreover, since $W_{\psi,\varphi}$ is invertible, it is clearly a Fredholm operator. It follows from \eqref{eq.3.18} and \eqref{eq.3.19} that
$$
0=\dim(\ker W_{\psi,\varphi})=N,\qquad
0=\codim(\Ran W_{\psi,\varphi})=K.
$$
Therefore, $\varphi$ is bijective and $G=T$. Thus condition \textup{(iii)} follows from Theorem~\ref{thm} condition \textup{(iv)}.

``$\Leftarrow$" Assume $W_{\psi,\varphi}$ satisfies conditions \textup{(i)--(iii)}. Then by Theorem~\ref{thm}, $W_{\psi,\varphi}$ is Fredholm. This yields that $W_{\psi,\varphi}$ has closed range.
Using \eqref{eq.3.18} and \eqref{eq.3.19} again, we see that
$$
\dim(\ker W_{\psi,\varphi})=\codim(\Ran W_{\psi,\varphi})=0.
$$
It follows that $W_{\psi,\varphi}$ is bijective.
Therefore, $W_{\psi,\varphi}$ is invertible on $\Lip$.
\end{proof}

The following example shows that there exists a Fredholm but non-invertible weighted composition operator on $\Lip$.

\begin{example}
Consider the tree $T=\mathbb{N}_0$, with $0$ being the root, where two vertices $n$ and $m$ are adjacent if and only if $|m-n|=1$.
Let
$$
\varphi(n)=n \ \ \text{for all} \ \ n\in\mathbb{N}_0
$$
and
\begin{align*}
\psi(n)=\begin{cases} 0, & \text{if } n=0,
\\ 1, & \text{if } n\geq 1.
\end{cases}
\end{align*}
By Theorem~\ref{thm}, $W_{\psi,\varphi}$ is Fredholm. However, since $Z_\psi\neq\varnothing$, it follows from Theorem~\ref{thm.4.1} that $W_{\psi,\varphi}$ is non-invertible on $\Lip$.
\end{example}

Next we characterize Fredholmness for composition operators on $\Lip$.
Let $\varphi: T\to T$. Put
$$
R=\varphi(T),\qquad F_w=\varphi^{-1}(w),\qquad N_w=|F_w|,
$$
$$
N=|T\setminus R|, \qquad K=\sum_{w\in R}(N_w-1)
$$
and
$$
G_\varphi=\{w\in R:F_w=\{s(w)\}\}.
$$
	
\begin{theorem}\label{thm.4.3}
Let $C_\varphi$ be a bounded composition operator on $\Lip$.
Then $C_\varphi$ is Fredholm if and only if the following conditions hold:
\begin{enumerate}[(i)]
\item $N<\infty$;
\item $|N_w|<\infty$ and $K<\infty$;
\item $\sup\{d(s(u),s(v)): u,v\in G_\varphi \ \text{with} \ \ u\sim v\}<\infty$.
\end{enumerate}
In this case,
$\ind C_\varphi=N-K$.
\end{theorem}

\begin{proof}
Using Theorem~\ref{thm} with $\psi\equiv1$, we have $Z_\psi=\varnothing$ and $G=G_\varphi$. If $u,v\in G_\varphi$ with $u\sim v$, put $x=s(u)$, $y=s(v)$ and $z=x\wedge y$. By the definition of $\Lambda_{\psi,\varphi}$, we obtain that
$$
\Lambda_{1,\varphi}(u,v)=d(z,x)+d(z,y)=d(x,y)=d(s(u),s(v)).
$$
Thus condition (iv) of Theorem~\ref{thm} is precisely condition (iii) of Theorem \ref{thm.4.3}.
This completes the proof.
\end{proof}

Finally we investigate the relationship between the Fredholmness of weighted composition operators and that of composition operators on $\Lip$. We begin with the simple lemma below, which will be used in our proof later.

\begin{lemma}\label{lem.4.4}
If $W_{\psi,\varphi}$ is Fredholm on $\Lip$, then
$$
\liminf_{|x|\to\infty}|\psi(x)|>0.
$$
\end{lemma}

\begin{proof}
By Theorem~\ref{thm}, we have
$$
M:=\sup_{\substack{u\sim v\\u,v\in G}}
\Lambda_{\psi,\varphi}(u,v)<\infty.
$$
Set $C=\max\{1,M\}$ and define
$$
B=\{o\}
\cup(T\setminus G)
\cup
\{v\in T:v\sim w\text{ for some }w\in T\setminus G\}.
$$
It follows that the set $B$ is finite. Let $u\in G\setminus B$. Then $u\neq o$ and $u^-\in G$. Put
$$
x=s(u),\qquad y=s(u^-),\qquad z=x\wedge y.
$$
Since $s$ is injective on $G$ and $u\neq u^-$, we have $x\neq y$. Hence at least one of $d(z,x)$ and $d(z,y)$ is positive.

If $d(z,x)\ge1$, then
$$
\frac1{|\psi(x)|}\le\frac{d(z,x)}{|\psi(x)|}\le\Lambda_{\psi,\varphi}(u,u^-)\le C
$$
and
$$
\left|\frac1{\psi(x)}-\frac1{\psi(y)}\right|
\le
\Lambda_{\psi,\varphi}(u,u^-)
\le C.
$$
It follows that
$$
\frac1{|\psi(y)|}
\le\frac1{|\psi(x)|}+\left|\frac1{\psi(x)}-\frac1{\psi(y)}\right|\le2C.
$$
If $d(z,y)\ge1$, with the similar argument, one can get the same conclusion. In particular, we have
\begin{equation}\label{eq.4.1}
|\psi(s(u))|\ge\frac1{2C},
\qquad u\in G\setminus B.
\end{equation}

If $x\notin s(G)$, then either $\psi(x)=0$ or $F_{\varphi(x)}$ contains more than one vertex. Therefore,
$$
T\setminus s(G)
\subseteq
Z_\psi
\cup
\bigcup_{\substack{w\in R\\N_w>1}}F_w.
$$
The right-hand side is finite by Theorem~\ref{thm}. Since $B$ is finite, we have $s(B\cap G)$ is finite, and hence the set
$$
T\setminus s(G\setminus B)
$$
is finite. By \eqref{eq.4.1}, $|\psi(x)|\ge\frac1{2C}$ for all but finitely many $x\in T$. This yields
$$
\liminf_{|x|\to\infty}|\psi(x)|>0.
$$
The proof is complete.
\end{proof}

\begin{theorem}
\label{thm.4.5}
Let $\varphi:T\to T$ and $\psi:T\to\mathbb C$. Then the following are equivalent:
\begin{enumerate}[(i)]
\item $W_{\psi,\varphi}$ is a Fredholm operator on $\Lip$;
\item both $M_\psi$ and $C_\varphi$ are Fredholm operators on $\Lip$.
\end{enumerate}
In this case,
$\ind M_\psi=0$ and $\ind W_{\psi,\varphi}=\ind C_\varphi$.
\end{theorem}

\begin{proof}
``$(ii)\Rightarrow(i)$" Obvious.

``$(i)\Rightarrow(ii)$" We use the same symbols as in Section 3. The proof is divided into several claims.

\begin{claim}
$C_\varphi$ is a Fredholm operator on $\Lip$.
\end{claim}

By Theorem~\ref{thm}, $N$, $K$, $|Z_\psi|$ are finite, and
$$
M:=\sup_{\substack{u\sim v\\u,v\in G}}
\Lambda_{\psi,\varphi}(u,v)<\infty.
$$
In particular, the set
$$
F=\bigcup_{\substack{w\in R\\N_w>1}}F_w
$$
is finite. If $x\in T\setminus F$, then
$F_{\varphi(x)}=\{x\}$. Choose a vertex $y$ adjacent to $x$.
Then $\varphi(y)\ne\varphi(x)$. For
$f=\chi_{\{\varphi(x)\}}$, we have $\|f\|_{\Lip}=1$ and
$$
(W_{\psi,\varphi}f)(x)=\psi(x),\qquad
(W_{\psi,\varphi}f)(y)=0.
$$
Consequently,
$$
|\psi(x)|\le \|W_{\psi,\varphi}f\|_{\Lip}
\le \|W_{\psi,\varphi}\|.
$$
Since $F$ is finite, it follows that
\begin{equation}\label{eq.4.2}
B:=\sup_{x\in T}|\psi(x)|<\infty.
\end{equation}
Moreover, we have $B>0$, because $Z_\psi$ is finite and $T$ is infinite.

By Lemma~\ref{lem.4.4}, there exist a finite set
$F_0\subset T$ and $m>0$ such that
$|\psi(x)|\ge m$ for all $x\in T\setminus F_0$.
Let $x\sim y$ and $c=\varphi(x)\wedge\varphi(y)$.
The functional
$$
f\longmapsto (W_{\psi,\varphi}f)(x)-(W_{\psi,\varphi}f)(y)
$$
has norm at most $\|W_{\psi,\varphi}\|$. Hence,
Lemma~\ref{lem.2.2} yields
$$
|\psi(x)|d(c,\varphi(x))
+|\psi(y)|d(c,\varphi(y))
\le \|W_{\psi,\varphi}\|.
$$
If $x,y\in T\setminus F_0$, this implies
$$
m\,d(\varphi(x),\varphi(y))\le \|W_{\psi,\varphi}\|.
$$
Note that $F_0$ is finite and $T$ is locally finite.
Then we have
$$
D_\varphi:=\sup_{x\sim y}d(\varphi(x),\varphi(y))<\infty.
$$
For any $f\in\Lip(T)$, it follows from
\eqref{eq.2.1} that
$$
\|C_\varphi f\|_{\Lip}
\le \max\{|\varphi(o)|+1,D_\varphi\}\|f\|_{\Lip}.
$$
This implies that $C_\varphi$ is bounded.

Furthermore, for any $u,v\in G$ with $u\sim v$, set
$x=s(u)$, $y=s(v)$, and $z=x\wedge y$.
By \eqref{eq.4.2}, we obtain that
$$
\Lambda_{\psi,\varphi}(u,v)
\ge \frac{d(z,x)+d(z,y)}{B}
=\frac{d(s(u),s(v))}{B}.
$$
It follows that $d(s(u),s(v))\le BM$.
Let $G_\varphi$ be the set given in
Theorem~\ref{thm.4.3}. Then $G\subseteq G_\varphi$ and
$$
G_\varphi\setminus G
=\{w\in G_\varphi:\psi(s(w))=0\}
\subseteq\varphi(Z_\psi),
$$
so this difference is finite.
Therefore
$$
\sup_{\substack{u\sim v\\u,v\in G_\varphi}}
d(s(u),s(v))<\infty.
$$
By Theorem~\ref{thm.4.3}, we obtain that $C_\varphi$ is
Fredholm with index $N-K$.

\begin{claim}
$M_\psi$ is a Fredholm operator on $\Lip$.
\end{claim}

Put $Y=\Ran C_\varphi$. Using
\eqref{eq.3.17} with weight identically equal to one, we see that
\begin{equation}\label{eq.4.3}
Y=\{h\in \Lip: h\text{ is constant on each } F_w\}.
\end{equation}
By Theorem~\ref{thm.4.3}, we have $N_w<\infty$ and
$K<\infty$. Thus, for each $F_w$ with $N_w>1$, we may fix
$b_w\in F_w$ and define
$$
Ph=\sum_{\substack{w\in R\\N_w>1}}
       \ \sum_{x\in F_w\setminus\{b_w\}}
       \bigl(h(x)-h(b_w)\bigr)\chi_{\{x\}},
\qquad h\in \Lip.
$$
It is clear that $P$ is finite rank and bounded. Moreover, on each $F_w$ with $N_w>1$,
$(I-P)h$ takes the constant value $h(b_w)$. Thus, by \eqref{eq.4.3}, we obtain
$(I-P)\Lip\subseteq Y$.

Since $C_\varphi$ has closed range, the bijection
$\Lip/\ker C_\varphi\to Y$ has bounded inverse. Hence there exists
$A>0$ such that
$$
d(f,\ker C_\varphi)\le A\|C_\varphi f\|_{\Lip},
\qquad f\in \Lip.
$$
For any $y=C_\varphi f\in Y$, we have
$M_\psi y=W_{\psi,\varphi}f\in\Lip$.
Since $\ker C_\varphi\subseteq\ker W_{\psi,\varphi}$, it follows that
$$
\begin{aligned}
\|M_\psi y\|_{\Lip}
&=\|W_{\psi,\varphi}f\|_{\Lip}\\
&\le \|W_{\psi,\varphi}\|
       d(f,\ker C_\varphi)\\
&\le A\|W_{\psi,\varphi}\|\|y\|_{\Lip}.
\end{aligned}
$$
Therefore, $M_\psi$ is bounded on $Y$.

Furthermore, we have
$$
M_\psi Ph=\sum_{\substack{w\in R\\N_w>1}}
\ \sum_{x\in F_w\setminus\{b_w\}}
\psi(x)\bigl(h(x)-h(b_w)\bigr)\chi_{\{x\}}.
$$
It follows that $M_\psi P$ is a bounded finite-rank operator on $\Lip$.
Consequently,
$$
\|M_\psi h\|_{\Lip}
\le \bigl(A\|W_{\psi,\varphi}\|\|I-P\|+\|M_\psi P\|\bigr)
\|h\|_{\Lip}
$$
for all $h\in\Lip$. This shows that $M_\psi$ is bounded on $\Lip$.

Finally, since both $M_\psi C_\varphi$ and $C_\varphi$ are Fredholm, $M_\psi$ is Fredholm on $\Lip$. By the additivity of the Fredholm index, we have
$$
\ind W_{\psi,\varphi}=\ind M_\psi+\ind C_\varphi.
$$
Therefore, using Theorem \ref{thm} and Theorem \ref{thm.4.3}, we see that
$$
\ind M_\psi=0,\qquad
\ind W_{\psi,\varphi}=\ind C_\varphi.
$$
\end{proof}

As an immediate consequence of Theorem~\ref{thm.4.5}, we obtain the following corollary.

\begin{corollary}
Let $\varphi:T\to T$ and $\psi:T\to\mathbb C$. Then the following are equivalent:
\begin{enumerate}[(i)]
\item $W_{\psi,\varphi}$ is an invertible operator on $\Lip$;
\item both $M_\psi$ and $C_\varphi$ are invertible operators on $\Lip$.
\end{enumerate}
\end{corollary}

{\bf Acknowledgement} This work is supported by the Natural Science Basic Research Program of Shaanxi Province(Program No. 2025JC-YBMS-024).

{\bf Author Contributions} All authors share credit for this work equally.

{\bf Data Availability} No datasets were generated or analysed during the current study.

{\bf Declarations}

{\bf Conflict of interest} The authors declare no competing interests.

\bibliographystyle{amsplain}

\end{document}